\documentclass[12pt]{article}
\usepackage{amsthm}
\usepackage{geometry}
\usepackage{graphicx}
\usepackage{comment}
\usepackage{amssymb}
\usepackage{amsmath}
\usepackage{fancybox}
\usepackage{epsfig}
\usepackage{palatino}
 \usepackage{setspace}
\usepackage{textgreek}

\usepackage{environ}
\NewEnviron{myequation1}{%
\begin{equation}
\scalebox{0.88}{$\BODY$}
\end{equation}}

\usepackage{tabularx, booktabs}

\begin{document}

\section*{Entropy-Variance Curves of Binary Sequences Generated by Random Substitutions of Constant Length}

\qquad

\noindent Juan Carlos Nuño and Francisco J. Mu\~noz \\
Department of Applied Mathematics\\
Universidad Politécnica de Madrid \\
Madrid 28040 (Spain) \\
juancarlos.nuno@upm.es

\qquad

\onehalfspacing

\section*{Abstract}
We study some properties of binary sequences generated by random substitutions of constant length. Specifically, assuming the alphabet $\{0,1\}$, we consider the following asymmetric substitution rule of length $k$: $0 \to \langle 0, 0, \ldots,0\rangle$ and $1 \to \langle Y_1, Y_2,  \ldots, Y_k \rangle$, where $Y_i$ is a Bernoulli random variable with parameter $p \in [0,1]$. We obtain by recurrence the discrete probability distribution of the stochastic variable that counts the number of ones in the sequence formed after a number $i$ of substitutions (iterations).  We derive its first two statistical moments, mean and variance, and the entropy of the generated sequences as a function of the substitution length $k$ for any successive iteration $i$, and characterize the values of $p$ where the maxima of these measures occur. Finally, we obtain the parametric curves entropy-variance for each iteration and substitution length. We find two regimes of dependence between these two variables that, to our knowledge, have not been previously described. Besides, it allows to compare sequences with the same entropy but different variance and vice versa.

\section{Introduction}

Binary sequences (strings or chains) appear naturally in physical systems for describing growing processes generated from an initial state. Examples are ubiquitous, from the spin up-down systems in physics to codification in information theory and its reflection in bio\-molecules as DNA and RNA (where two sets of nucleotides exist: purines and pyrimidines).  They are also common objects in mathematics, representing, for example, sets of words formed from letters of an alphabet or symbolic representations of dynamical systems, where mapping can be defined to transform real trajectories into two state series \cite{Fogg}.  Binary sequence can also codified fractal sets obtained by an iterative process, as occurs in the classical Cantor set \cite{Falconer,Dekking}. 

In order to study the properties of systems with fractal geometry, Mandelbrot defined some random sets that are generated recursively from a given initial set as a percolation process \cite{Mandelbrot, Dekking}. These sets are formed by successive application of a defined set of rules that, either divide the initial set and discard some subsets in each partition, or enlarge the initial set by substituting the existing subsets by multiple replicas of themselves. These sets can be embedded in 1D, 2D and 3D spaces, conditioning some of their characteristics, in particular, percolation.

In this paper, we study the properties of the sequences formed by performing random substitutions of constant length (see, for instace, \cite{Barbe}).  Specifically, we consider a binary alphabet, that we denotate as $\{0,1\}$, and substitutions of length $k=2,3, \ldots $. As in the classical percolation process, we assume that, once a void site is formed, i.e., a 0 occurs at this site, subsequent substitutions of it yield $k$ zeros, always. On the other hand, filled sites, i.e., 1 at this site, are substituted by random words of $k$ letters with a (uniform) probability of inserting a 1 defined by a parameter $p$. In this substitution, the probability of inserting a 0 is $q=1-p$. 

The process defined in this way gives rise to sequences of increasing size. The number of ones (filled sites) in the sequence at iteration $i$ is a stochastic variable that depends on the probability $p$ and the length $k$. We study the probabilistic properties of these sequences. Specifically, we derive the probability distribution of the number of ones at iteration $i$. We also obtain the expected values of $X$ and $X^2$, the corresponding variance $VAR$, and obtain some properties of them.  As an alternative measure of uncertainty, we calculate the mean entropy of these sequences and compare it with their variance \cite{Ebrahimi, Mukherjee}. Finally, we find the parametric curves that relate both measures for each iteration and substitution length
. These curves present two regimes: for low values of $p$ the dependence entropy-variance is convex (concave down), whereas for large values of $p$ this dependence is practically lost. Therefore, for each iteration and substitution length, we can generate sequences with the same entropy but with different variance and vice versa.

\section{Sequence Generation}\label{sec2}

\begin{figure}[h]
\centering
\includegraphics[width=1\textwidth]{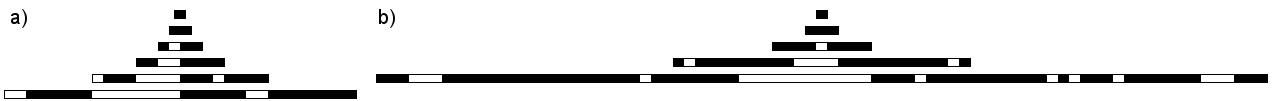}
\caption{An example of schematic representation of the sequence generation (from top) for the cases ({\bf a}) $k=2$ and ({\bf b}) $k=3$, after six and five substitutions, respectively. For both cases, $p=0.9$. Black and white represent 1 and 0, respectively.}
\label{Figgen}
\end{figure}

The generation of the sequences studied in this work is a one-dimensional example of the classical well-known Mandelbrot's percolation process \cite{Mandelbrot}. Contrary to the classical formulation, where an initial segment of length $L$ is subdivided into $k$ subsegments of length $L/k$ and, with probability $p$, some of them are chosen to be further divided, the model we study here considers an initial digit, that we denotate by $1$, that is substituted by a word of length $k$ formed by digits  $\{0,1\}$.  As in the classical model, we assume that digits denotated by 0 are substituted by a word of $k$ zeros (see Figure \ref{Figgen}) and sites occupied by ones are substituted by a word of $k$ letters with a independent probability $p$. This procedure is applied successively to generate a binary sequence of finite length, after $i$ iterations or, in the limit, infinite length.  When $k=3$, this procedure can be viewed as a probabilistic version of the Cantor set for $0<p<1$ \cite{Dekking}. 

These procedures are called substitutions, as they are obtained by replacing a digit by a word of digits \cite{Morales,Fogg, Wolfram}. In this context, we study the following map that applies at each random substitution of length $k$:
\begin{eqnarray}\label{Grule}
0 & \to &  \langle 0 \, 0 \, \ldots \, 0 \rangle\\ \nonumber
1 & \to &  \langle Y_1 \, Y_2 \, \ldots \, Y_k  \rangle
\end{eqnarray}
where $Y_j$ is a Bernoulli random variable with parameter $p \in [0,1]$, i.e., $p$ means the (uniform) probability of inserting 1 in each position and independently for all $j$. We note that an alternative description of this substitution process can be defined by using the Kronecker product \cite{Glebeta, Voevudko,Shallit,Xue}.

The deterministic version of this procedure brings about well known sequences.  For the classical Cantor set, this map reduces to: 
\begin{eqnarray}
0 & \to &  \langle 0 \, 0 \, 0 \rangle \\ \nonumber
1 & \to &  \langle 1 \, 0 \, 1 \rangle
\end{eqnarray}
which, starting from a unique 1, yield the sequences:
\begin{eqnarray}
s_1& =&  (1,0,1) \\ \nonumber
s_2& =&  (1,0,1,0,0,0,1,0,1) \\
s_3 & = & (1,0,1,0,0,0,1,0,1,0,0,0,0,0,0,0,0,0,1,0,1,0,0,0,1,0,1) \nonumber
\end{eqnarray}
whose limit of the iterative process is the infinite sequence that corresponds to the classical Cantor set. For this particular case,  the length of the $i$-th vector is $L_i=3^i$. In general, for substitutions of constant length $k$,  the $i$th-iterated vector has length $L_i=k^i$. 

{There are other well-known sequences generated iteratively by substitutions of different length. For instance, the Fibonacci sequence is generated applying the following map:
\cite{Schroeder, Baake}:}
\begin{equation}
0 \to \langle 1 \rangle \hspace*{1cm} 1 \to \langle 1 \,0 \rangle
\end{equation}
Starting from 0, after six iterations we obtain the binary vector: $(1,0,1,1,0,1,0,1)$. 
Here, the number of 1 of the $i$th-sequence represents the fertile population of rabbits. Thus, the number of mature pairs of rabbits after six generations is 5. Notice that, in principle, the position in the sequence has not a spatial meaning. It is also worthy to remark that this sequence can be generated by different maps \cite{Schroeder}.

Another self-similar sequence generated by a recursive process is the so-called Morse-Thue (MT), which is obtained by the following map:
\begin{equation}
0  \to  \langle 0 \, 1 \rangle \hspace*{1cm}  1 \to  \langle 1\, 0 \rangle
\end{equation}
starting from 0 (although it could be also initiated by 1) \cite{Schroeder,Baake}. Also for this example, there is an alternative way of generating this sequence, simply by take the sequence of the preceding step and add the complement. Besides, this sequence is aperiodic, although completely deterministic. On the contrary, it can be proven that there are short and long range correlations \cite{Schroeder}.

\section{Main Statistical Properties}

In this section, we study the statistical properties of the sequences generated by applying the rule (\ref{Grule}) recurrently to the starting sequence ($ 1$). A classical issue to be considered is the probability distribution of each of the letters of the alphabet. In particular, in the case of binary sequences, this problem is reduced to know the probability distribution of, for instance, the number of ones of any population of sequences generated from a set of random rules, e.g. random substitutions.

As a starting point, the classical Bernoulli process could be illustrative. This process can be also viewed as a random substitution of constant length, $k$, where both letters, e.g., 0 and 1, are identically substituted according with the rule:
\begin{equation}
0 \to \langle Z_1, \, Z_2 \, \ldots, \, Z_k  \rangle; \hspace*{1cm} 1 \to \langle Y_1, \, Y_2 \, \ldots, \, Y_k  \rangle
\end{equation}
being $Z_j$ and $Y_j$ Bernoulli processes with parameter $p \in [0,1]$. In each iteration, $i$, the length of the sequence is $n = k^i$, and  the probability distribution of the number of ones, $X$, is given by the binomial distribution:
\begin{equation}\label{Berdistri}
 P(X=x) = \binom{n}{x} \, p^x (1-p)^{n-x}
\end{equation}
By applying this distribution function, the expected values of the number of ones in a sequence of length $n$, $X_n$, and its square, respectively $E ^B(X_n)$ and $E^B(X_n^2)$, can be calculated:
\begin{eqnarray}\label{expecB}
E^B ( X_n)  & = & n \, p  \\
E^B ( X_n^2) & =& n \,p \, (1 + (n-1) \, p)
\end{eqnarray}
for $n=1, 2, 3, \ldots$. Using these expressions, the variance of $X$ reads:
\begin{equation}
VAR^B (X_n)  =E^B (X_n^2)-(E^B (X_n ))^2 = n \, p \, (1-p)
\end{equation}
Note that for the binomial distribution, the ratio variance-mean is equal to $1-p$ for any $n$.

Similarly, we can compute these moments for the sequences generated from random substitutions defined in the previous section. In this case, a recurrent formula 
provides this distribution for any iteration $i$ and for any value of the substitution length $k$. 

\begin{figure}[h]
\includegraphics[width=0.5\textwidth]{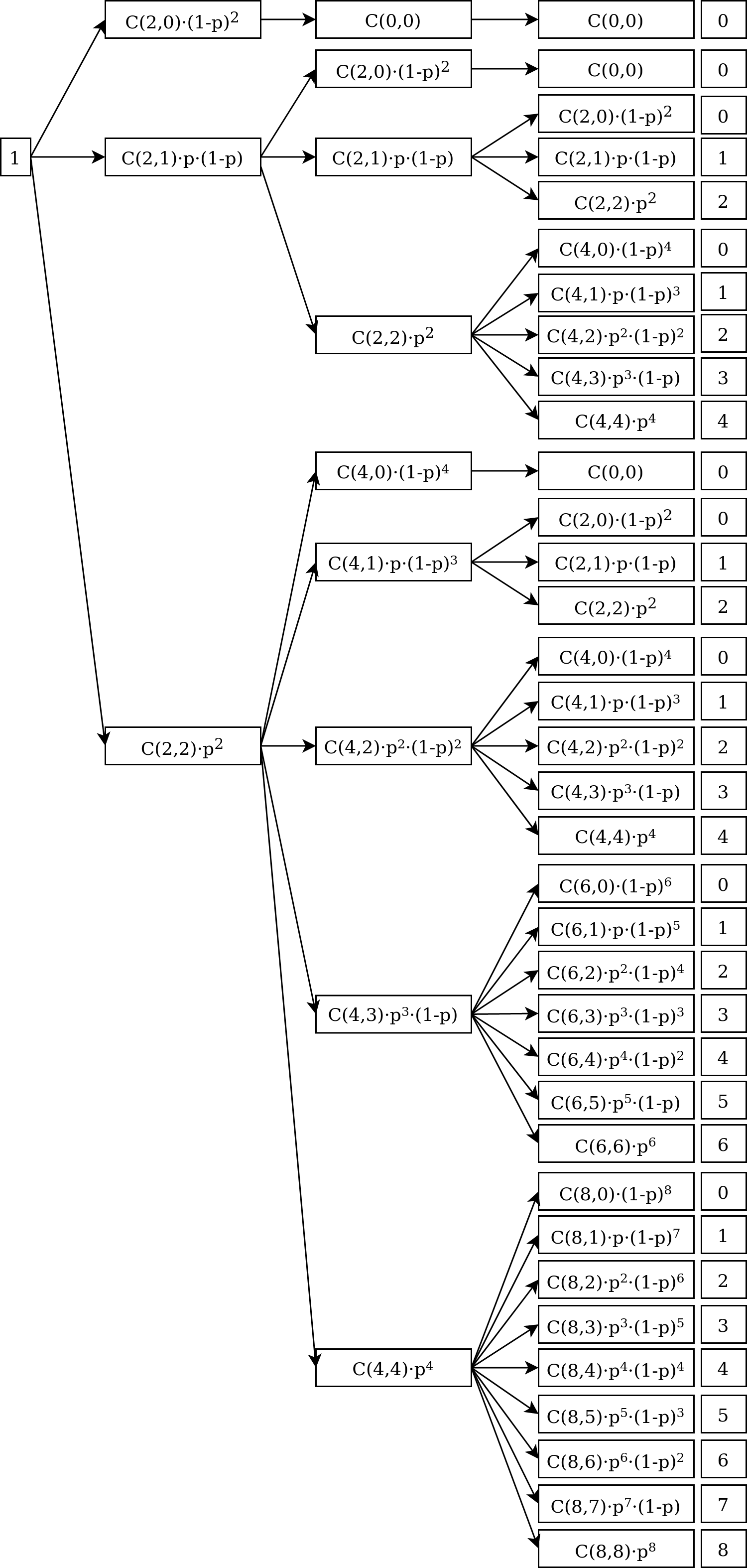}
\caption{A tree representation of the probabilities of obtaining a number of ones after $i=3$ iterations using rule (\ref{Grule}). As usual, $C(i,j)$ means the combinatorial numbers $C(i,j)=\binom {i}{j}$. The rightmost column represents the number of ones for each branch.} 
\label{Ptree}
\end{figure}

Let us present first the results for a substitution of length $k=2$. As stated before, we start the generation of the binary sequences from a sequence with a unique one. Thus, the probability of having this initial sequence is assumed to be 1:
\begin{equation}
P(X_0=1) = 1
\end{equation}
In the first step, a sequence of length two is formed whose distribution of ones depends on the probability $p$ (see Figure \ref{Ptree}). To obtain the the distribution probability of ones in this first step, $P({\bf X}_1)$, we multiply from left $P(X_0)$ by the substitution matrix:
\begin{equation}\label{M31}
M_{3\times1} = 
\begin{pmatrix} 
\binom{2}{0} \, (1-p)^2 \\ 
\binom{2}{1} \, p\,(1-p) \\ 
\binom{2}{2} p^2 
\end{pmatrix}  
\end{equation}
givin rise to the vector:
\begin{equation}\label{1ite}
{\bf P}({\bf X}_1) = 
\begin{pmatrix} P_1(X=0) \\ P_1(X=1) \\ P_1(X=2) \end{pmatrix} = 
 M_{3\times1} \, P(X_0)  =  M_{3\times1}
\end{equation}
where ${\bf X}_1 = (0,1,2)$.

In the next substitution, a four digits sequence is formed and the probability distribution of ones can be computed from the previous one as follows:
\begin{myequation1}\label{2ite}
{\bf P}({\bf X}_2) = \begin{pmatrix} P (X=0) \\ P (X=1) \\ P(X=2) \\ P(X=3) \\ P(X=4) \end{pmatrix} = 
 \begin{pmatrix} 
 1 & \binom{2}{0} (1-p)^2 & \binom{4}{0} \, (1-p)^4 \\
0 & \binom{2}{1} p\,(1-p) & \binom{4}{1} \, p\, (1-p)^3 \\
0 & \binom{2}{2} p^2 & \binom{4}{2} \, p^2\, (1-p)^2 \\
0 & 0  & \binom{4}{3} \, p^3\, (1-p) \\
0 & 0  & \binom{4}{4} \, p^4 
\end{pmatrix}
\begin{pmatrix}
P(X=0) \\ P(X=1) \\ P(X=2) 
\end{pmatrix}=
M_{5\times3}\, {\bf P}({\bf X}_1)
\end{myequation1}
where ${\bf X}_2 = (0,1,2,3,4)$. 
The third iteration confirms the formula (see Figure \ref{Ptree}):
\begin{equation}
{\bf P }({\bf X}_3) = M_{9\times5}\,{\bf P}({\bf X}_2)
\end{equation}
with 
\begin{equation}
{\bf P} (\bf{X}_3) = 
\begin{pmatrix} P(X=0) \\ P(X=1) \\ P(X=2) \\ P(X=3) \\ P(X=4) \\ P(X=5) \\ P(X=6) \\ P(X=7) \\ P (X=8) \end{pmatrix}
\end{equation}
where ${\bf X}_3 = (0,1,2,3,4,5,6,7,8)$, and 
\begin{equation}
M_{9\times5} = 
\begin{pmatrix}
 1 & \binom{2}{0} (1-p)^2 & \binom{4}{0} \, (1-p)^4 &  \binom{6}{0} (1-p)^6 & \binom{8}{0} \, (1-p)^8 \\
 0 & \binom{2}{1} p\, (1-p) & \binom{4}{1} \, p\, (1-p)^3 &  \binom{6}{1} p\,(1-p)^5 & \binom{8}{1} \, p\,(1-p)^7 \\
 0 & \binom{2}{2} p^2 & \binom{4}{2} \, p^2\, (1-p)^2 &  \binom{6}{2} p^2\,(1-p)^4 & \binom{8}{2} \, p^2 \,(1-p)^6 \\
 0 & 0 & \binom{4}{3} \, p^3\, (1-p) &  \binom{6}{3} p^3\,(1-p)^3 & \binom{8}{3} \, p^3 \,(1-p)^5 \\
 0 & 0 & \binom{4}{4} \, p^4  &  \binom{6}{4} p^4\,(1-p)^2 & \binom{8}{4} \, p^4 \,(1-p)^4 \\
 0 & 0 & 0  &  \binom{6}{5} p^5\,(1-p) & \binom{8}{5} \, p^5 \,(1-p)^3 \\
 0 & 0 & 0  &  \binom{6}{6} p^6  & \binom{8}{6} \, p^6 \,(1-p)^2 \\
 0 & 0 & 0  &  0 & \binom{8}{7} \, p^7 \,(1-p) \\
 0 & 0 & 0  &  0 & \binom{8}{8} \, p^8 
\end{pmatrix}
\end{equation}

The general expression for the recurrent formula for the probability distribution of the number of ones in the $i$-th iteration is:
\begin{equation}\label{nite}
{\bf P} ({\bf X}_i) = M_{(2^i+1)\times(2^{i-1}+1)}\,{\bf P} ({\bf X}_{i-1})
\end{equation}
being the vector ${\bf X}_i=(0,1,2,\ldots,2^i)$ and the general substitution matrix, $M_{(2^i+1)\times(2^{i-1}+1)}$, whose coefficients are given by: 
\begin{equation}
m_{ab}^i = \left \{ \begin{matrix} \binom{2\,(b-1)}{a-1} \, p^{a-1} \, (1-p)^{2\,(b-1) - (a-1)} & \mbox{if} &  a \leq 2\,b -1  \\
0 &   & \mbox{otherwise} 
\end{matrix}
\right.  
\end{equation}
for $i=2,3,\ldots$ and $1 \leq a \leq 2^i +1$ and $1 \leq b \leq 2^{i-1} +1$. 
For $i =1$, the coefficients of matrix $M_{3\times1}$ (\ref{M31}) are given by:
\begin{equation}
m_{a1}^1 = \binom{2}{a-1} \, p^{a-1} \, (1-p)^{2 - (a-1)} \hspace{12pt} \mbox{for} \hspace{8pt}  1\leq a \leq 3  
\end{equation}
By recurrence, Formula \eqref{nite} can be converted into a matrix product chain: 
\begin{equation}\label{Mexp}
{\bf P} ({\bf X}_i) =  M_{(2^i+1)\times(2^{i-1}+1)}\,  M_{(2^{i-1}+1)\times(2^{i-2}+1)} \ldots  M_{9\times5}\, M_{5\times3}\,{\bf P} ({\bf X}_1)
\end{equation}

\begin{figure}[h]
\includegraphics[width=0.9\textwidth]{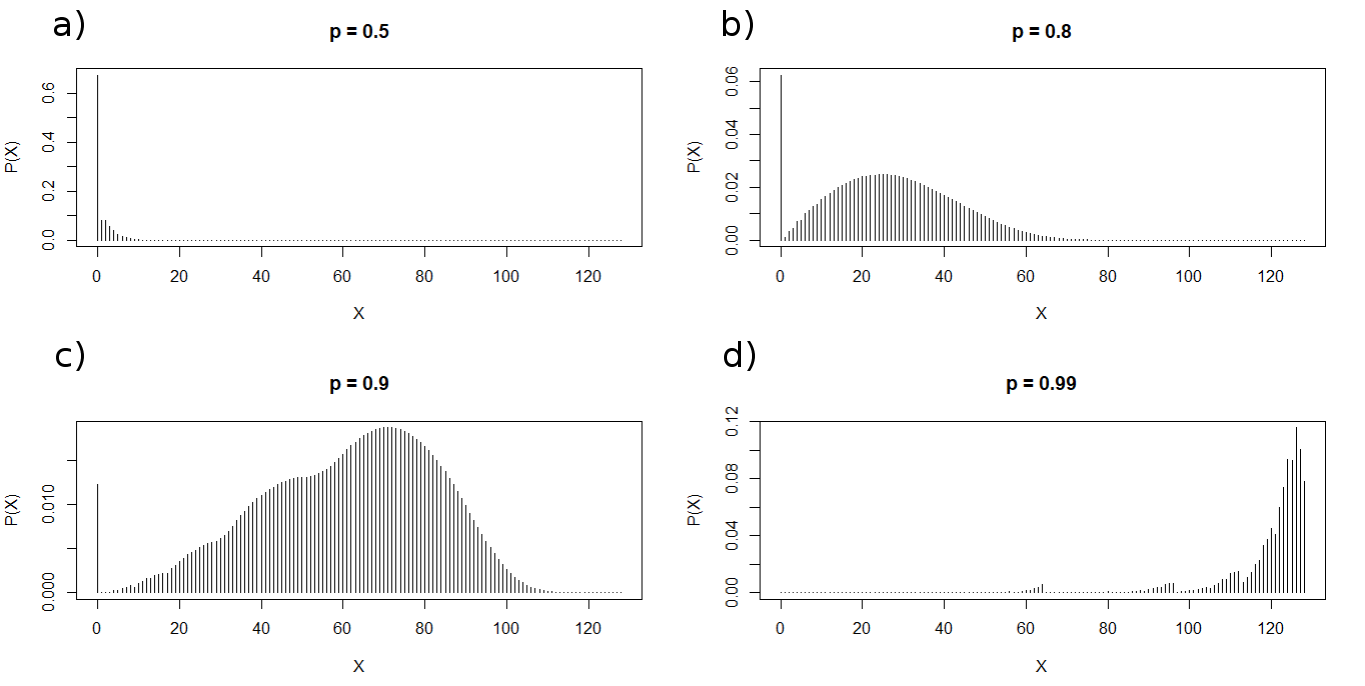}
\caption{Discrete probability distributions of the number of ones for different values of $p$ for the seventh iteration ($i=7$). Concretely, $p=0.5$ ({\bf a}) , $p=0.8$ ({\bf b}) , $p=0.9$ ({\bf c}) and $p=0.99$ ({\bf d}) . As it can be seen, the distributions for low value of $p$ exhibit a peak at $X=0$ that disappears when $p$ tends to 1. }
\label{Probdistri}
\end{figure}

\begin{figure}[h]
\includegraphics[width=0.9\textwidth]{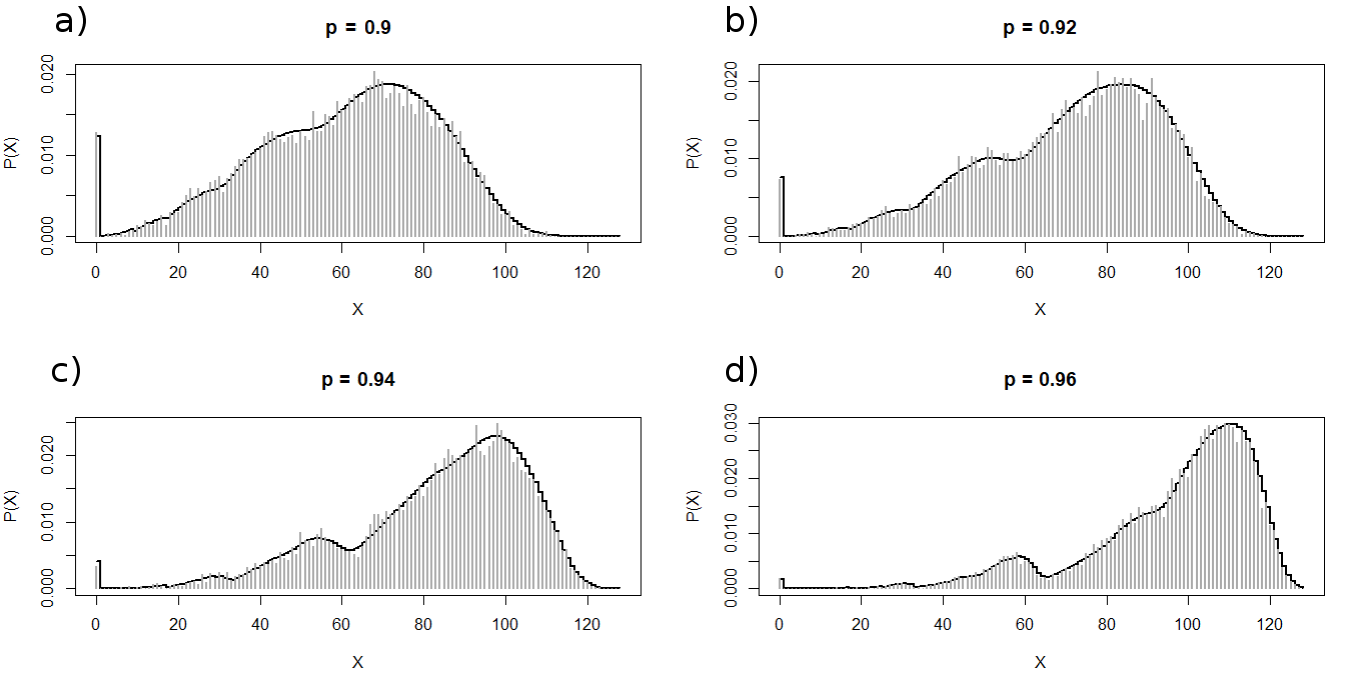}
\caption{Analytical (black continuous curve) and the corresponding numerical distributions for $p=0.9$ ({\bf a}), $p=0.92$ ({\bf b}), $p=0.94$ ({\bf c}) and $p=0.96$ ({\bf d}). The numerical distributions are histograms (gray bars) obtained from 1000 realizations. Note the multimodal character of the distributions; for large values of $p$, in addition to the peak at $X=0$, other peaks appear at $X \approx 18, 30, 54$ and $85$. The heights of these maxima vary with $p$ and new peaks are visible as $p$ increases (although, in the limit $p=1$, only one remains).}
\label{Probdistri2}
\end{figure}

Figure \ref{Probdistri} depicts the probability distributions for different values of $p$ for the iteration $i=7$. Figure \ref{Probdistri2} shows both the analytical (black continuous curves) and the corresponding numerical distributions obtained from 1000 realizations. As it can be seen, the distributions exhibit a peak at $X=0$ whose height decreases as $p$ increases. Besides this peak, the distributions also present other maxima whose heights depends on $p$, being smaller for lower values of $p$. This multimodality could be a consequence of the generation rule that fosters the emergence of clusters of scalable sizes. As a matter of fact, the fractal properties of these sequences are the main subject of a forthcoming paper.

The probability of generating a null sequence, i.e., formed completely by 0,  after $i$ substitutions, {$P (X_i=0)$}
, corresponds to the first component of {${\bf P} ({\bf X}_i)$}. The following recurrence formula can be obtained \cite{Chayes}:
\begin{equation}\label{Prob0}
P (X_{i+1}=0) = \left(P (X_i=0)\, p + (1-p) \right)^2
\end{equation}
for $i=1,2,\ldots$.  Note that  $p\,P(X_i \neq 0)$ represents the probability of yielding a non-null sequence after a substitution that includes at least one 1 in the previous non-null sequence and then, $1 - p \, P(X_i \neq 0)$ represents the probability of yielding a null sequence after a substitution of length 1. This expression can be rewritten as:
\begin{equation}
1 - p \, P(X_i \neq 0)  = 1 - p (1 - P(X_i = 0) ) = 1- p + p \, P(X_i = 0)
\end{equation}
Because the substitution has length $k=2$, and it is formed by independent digits, this expression must be multiplied by itself resulting in Equation \eqref{Prob0}.

The first term in this sequence is: 
\begin{equation}
P (X_1=0) = (1-p)^2 
\end{equation}
which yields the second one:
\begin{equation}
P (X_2=0) = \left((1-p)^2\,p + (1-p)\right)^2 
\end{equation}

It is worthy to remark that the expansions of $P_i(X_i=0)$ as polynomials of $p$ and $q=1-p$ have the {\it triangle by rows} (oeis-A202019) \cite{oeis} as coefficients. 

\begin{figure}[h]
\includegraphics[width=0.5\textwidth]{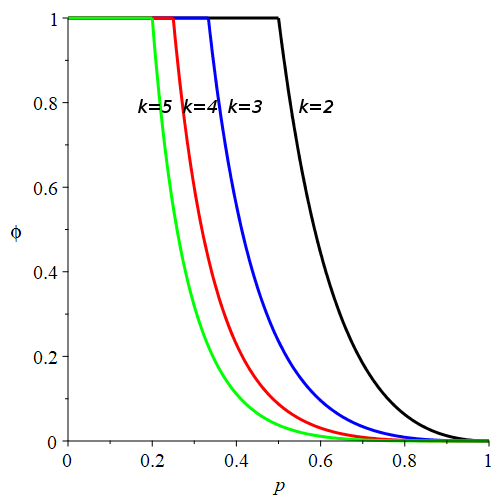}
\caption{Probability distribution of having a null sequence at the limit of the generating process, $\phi$, for $k=2$ (black), $k=3$ (blue), $k=4$ (red) and $k=5$ (green), as a function of $p$.}
\label{figphi}
\end{figure}

Assuming that this sequence of probability distributions converges to $\phi$, the following equation holds:
\begin{equation}
\phi = \left(\phi \, p + (1-p)\right)^2
\end{equation}
whose lowest solution: 
\begin{equation}
\phi(p) = \min\{1,\left(\frac{1-p}{p}\right)^2\}
\end{equation}
depends on $p$ as depicted in Figure \ref{figphi}. Note that,  for $p$ below $p_c=\frac{1}{2}$ the probability of generating a null infinite sequence is one, i.e.,  $\phi(p) = 1$ for $p<p_c$.

Equation \eqref{Prob0} can be generalized straightforwardly to any substitution length $k$ using the same demonstration as above (see \cite{Chayes}). Now, the probability of having a null sequence at iteration $i+1$ is obtained after multiplying $k$ times the factor $1 - p \, P(X_i \neq 0)$ :
\begin{equation}\label{Prob0k}
P (X_{i+1}=0) = \left(1- p \, P(X_i \neq 0) \right)^k = \left(P (X_i=0)\, p + (1-p) \right)^k
\end{equation}

As before, in the limit this equation converges to:
\begin{equation}
\phi = \left(\phi \, p + (1-p)\right)^k
\end{equation}
that also exhibits a critical value of $p$ :
\begin{equation}
p_c(k) = \frac{1}{k}
\end{equation}
As it can be observed in Figure \ref{figphi}, as the substitution length $k$ increases the probability of generating a null sequence decreases for $p>p_c(k)$; also does the range where this probability is one, i.e., $p_c(k) \to 0$ as $k \to \infty$.

\subsection*{Expected Values and Variance}

Having already the probability distributions of ones, it is straightforward to compute its expected value after $i$ iterations:
\begin{equation}
E ({\bf X}) = \sum_{j=0}^{2^i} j \, P_i (X = j) 
\end{equation}
This equation can be written in matrix form:
\begin{equation}
E ({\bf X}_i) = {\bf X}_i \, {\bf P} ({\bf X}_{i}) 
\end{equation}
Now, by applying Formula \eqref{nite}
\begin{equation}
E ({\bf X}_i) = {\bf X}_i \, M_{(2^i+1)\times(2^{i-1}+1)}\,{\bf P} ({\bf X}_{i-1}) 
\end{equation}
It turns out that
\begin{equation}
{\bf X}_i \, M_{(2^i+1)\times(2^{i-1}+1)} = { E^B}({\bf X}_i)
\end{equation}
where by ${ E^B} ({\bf X}_i)$ we refer to a vector whose entries are the expectation of a Binomial random variable $Y_j \sim Bin(2\,j,p), \, j=0,1,2,...,2^{i-1}$ for successive iterations. Thus:
\begin{equation}
{ E^B}({\bf X}_i) = 2 \, {\bf X}_i = (0, 2\,p, 4\, p , \ldots, 2^i \, p)
\end{equation}
Then, the following recurrence formula holds:
\begin{equation}
E ({\bf X}_{i+1}) = 2\,p \, E ({\bf X}_i)
\end{equation}
for $i=0,1,2,\ldots$ and $E_1=2\,p$.  This recurrence yields the expression:
\begin{equation}\label{meanmu}
E ({\bf X}_i) = \mu_i(p,2) = (2\,p)^i
\end{equation}
for $i=0,1,2,\ldots$.

This formula can be straightforwardly generalized to any length $k$: 
\begin{equation}\label{meanmuk}
E ({\bf X}_i) = \mu_i(p,k) = (k\,p)^i \hspace*{5mm} i =0,1,2, \ldots
\end{equation}

The mean of the distribution of 0 in the sequence formed after $i$ substitutions of length $k$ is:
\begin{equation}
\mu_i^0(p,k) = k^i - (kp)^i = k^i (1-p^i)
\end{equation}
that yields a ratio $\#1/\#0$ at the $i$th-iteration:
\begin{equation}
\mu_i^{10}(p) = \frac{(kp)^i}{ k^i \, (1-p^i)} = \frac{p^i}{1 - p^i}
\end{equation}
This ratio tends to 0 as $i \to \infty$ for all $0 \leq p < 1$ and for all $k$. Furthermote, the limit of $\mu_i^{10}(p)$ as $p \to 1$ is infinite for all values of $i=1,2,\ldots$. Note that this ratio is one for the value of $p_{\mu=1} = 2^{- \frac{1}{i}}$, independently of $k$. 

Similarly, it can be calculated the expected value of ${\bf X}^2$, $E({\bf X}^2)$. First, for the case $k=2$ we obtain:
\begin{equation}
E ({\bf X}_i^2) = \sum_{j=0}^{i} j^2 \, P (X = j) 
\end{equation}
For the second iteration, it is straightforward to calculate this value from this definition:
\begin{equation}
E ({\bf X}_1^2)=(0,1,4) \begin{pmatrix} (1-p)^2 \\ 2\, p\,(1-p) \\ p^2 \end{pmatrix} = 2\,p \, (1+p)
\end{equation}
Using the recurrence of the probability distributions for successive iterations, the expected value of the next iteration is also computed:
\begin{equation}
E ({\bf X}_2^2) = (0,1,4,9,16)  \,  \begin{pmatrix} 
 1 & \binom{2}{0} (1-p)^2 & \binom{4}{0} \, (1-p)^4 \\
0 & \binom{2}{1} p\,(1-p) & \binom{4}{1} \, p\, (1-p)^3 \\
0 & \binom{2}{2} p^2 & \binom{4}{2} \, p^2\, (1-p)^2 \\
0 & 0  & \binom{4}{3} \, p^3\, (1-p) \\
0 & 0  & \binom{4}{4} \, p^4 
\end{pmatrix}
\begin{pmatrix} (1-p)^2 \\ 2\, p\,(1-p) \\ p^2 \end{pmatrix} 
\end{equation}
that simplifies to: 
\begin{equation}
E ({\bf X}_2^2)= (0, 2\,p\,(1+p), 4\,p\, (1+3p)) \begin{pmatrix} (1-p)^2 \\ 2\, p\,(1-p) \\ p^2 \end{pmatrix} = 4 \, p^2 (2\,p^2 + p + 1)
\end{equation}
Here we have used the expected value of $X^2$ for the binomial distribution for the number of ones in a sequence of length $n$.

In order to calculate the expected values of $X^2$ for successive iteration we find the following recursive formula:
\begin{equation}\label{REc2}
E ({\bf X}_i^2) = 2\,p\,(1-p) \, E ({\bf X}_{i-1})  + (2\, p)^2 \, E\, ({\bf X}_{i-1}^2) 
\end{equation}
that can be derived using Equation \eqref{expecB} for the expected value of $X^2$. Effectively, 
\begin{equation}
E({\bf X}_i^2) = {\bf X}_i^2 \, M_{(2^i+1)\times(2^{i-1}+1)} \, {\bf P}({\bf X}_{i-1}) =  E^B({\bf X}_i^2) \, {\bf P} ({\bf X}_{i-1}) 
\end{equation}
where 
\begin{equation}
{\bf X}_i^2 = \left(0, 1, 2^2 , 3^2, \ldots,  (2^i)^2  \right)
\end{equation}
As before, $E^B({\bf X}_i^2)$ represents a vector with entries $E(Y_j^2)$ where $Y_j \sim Bin(k\,j,p)$ \newline for \mbox{$j=0,1,2,...,2^{i-1}$} and, consequently, is given by:
\begin{equation}
E^B({\bf X}_i^2) = \left(0, 2\,p \, (1+p), 2^2 \, p \, (1 + 3 \, p), \ldots,  2^i \, p \, (2^{i}-1)\, p) \right) 
\end{equation}
Then, since this vector can be split into two vectors:
\begin{equation}
E^B({\bf X}_i^2) = \left(0, 2\,p, 2^2 \, p , \ldots,  2^i \, p  \right) + \left(0, 2\,p^2 , 2^2 \, p^2 , \ldots,  2^i \, p^2 \right) 
\end{equation}
Equation \eqref{REc2} holds.

It can be proven that the expected value of ${\bf X}^2$ for the $i$-th iterations depends on $p$ as follows:
\begin{equation}
E ({\bf X}_i^2) =  2^{i-1} \, p^i  \, \left( 1 + \sum_{j=0}^i  (2\,p)^j \right)
\end{equation}
for $n=1,2,3,\ldots$.  If we replace the geometric sum:
\begin{equation}
E ({\bf X}_i^2) =\frac{1}{2} \, (2\,p)^i \, \left(1 + \frac{1-\, (2\,p)^{i+1}}{1-2\,p}\right)
\end{equation}

This expression can be applied to calculate the variance of the distribution of the number of ones in a sequence generated after $i$ iterations, i.e., $VAR ({\bf X}_i) = E({\bf X}_i^2) - (E({\bf X}_i))^2$:
\begin{equation}
VAR ({\bf X}_i) =\frac{1}{2} \, (2\,p)^i \, \left(1 + \frac{1-\, (2\,p)^{i+1}}{1-2\,p}\right) - (2\,p)^{2\,i} 
\end{equation}
that simplifies to:
\begin{equation}\label{GVAR}
VAR ({\bf X}_i) = \frac{1-p}{1-2\,p} (2\,p)^i (1-(2\,p)^i)
\end{equation}
This expression can be further compacted if we use the expression of the mean (\ref{meanmu}):
\begin{equation}\label{sigmu}
VAR ({\bf X}_i) = \frac{1-p}{1-2\,p} \, \mu_i(p) (1-\mu_i(p))
\end{equation}

These computations can be equally performed for any length $k =3,4,\ldots$. A general formulation for any value of $k$ yields the following functions of $p$ for the expected values of $X$ and $X^2$ at the $i$-th iteration:
\begin{eqnarray}
E ({\bf X}_i) (p,k)& = & (k\,p)^i \\
E ({\bf X}_i^2)(p,k) & =  & (k\,  p)^i \, \left(1 + (k-1) \, p \, \frac{1- (k\,p)^i}{1-k \,p}\right)
\end{eqnarray}
Then, the variance is straightforward calculated and yields:
\begin{equation}
VAR ({\bf X}_i) = \frac{1-p}{1-k\,p} (k\,p)^i (1-(k\,p)^i)
\end{equation}

The limits as $i$ tends to infinity of the mean $\mu$ and the variance $VAR$ depend on the value of $p$. For $0< p <\frac{1}{k}$, these limits equal 0, whereas for $\frac{1}{k} < p < 1$, the limits are~infinite.

\begin{figure}[h]
\includegraphics[width=0.5\textwidth]{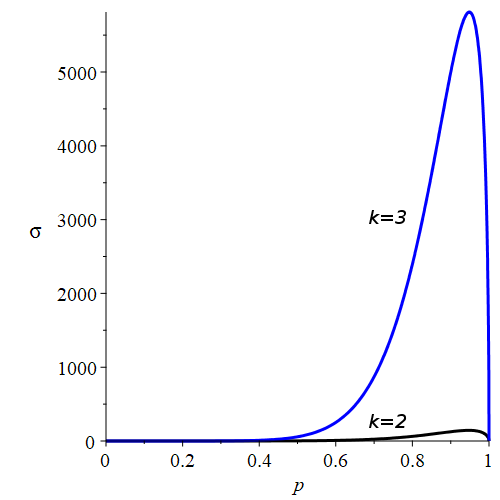}
\caption{Standard deviation, $\sigma= \sqrt{VAR({\bf X}_i)}$ for the substitution lengths $k=2$ (black) and $k=3$ (blue) after $i=10$ substitutions as a function of $p$. Note the different scale at the $Y$-axis, likely as a consequence of the different sequence size $2^{10}$ vs. $3^{10}$.}
\label{varphi}
\end{figure}

Contrary to the Bernoulli sequences, the variance reaches a maximum for an intermediate value of $p$, that depends on both $k$ and $i$ (see Figure \ref{varphi}).

An interesting index that could provide information about the stochasticity of the generation process is the variance-mean ratio, also known as the dispersion index, $D$. For iteration $i$, it can be written as:
\begin{equation}
D_i(p,k) =\frac{VAR ({\bf X}_i)}{E ({\bf X}_i)} = \frac{1-p}{1-k\,p} (1-(k\,p)^i)
\end{equation}
By comparison with the Poisson distribution, that possesses  a value of the dispersion equal to 1, it is said that a process that has a  $D$-value lower than 1 is ``under-dispersed''; on the contrary, having a $D$-value larger than 1 reflects an ``over-dispersed'' stochasticity. Since $D_i(0,k) = 1$ and the derivative of $D_i(p,k)$ is positive at $p=0$ for all $i>1$ and furthermore, $D_i(1,k) = 0$, then there exist a value of $p (k,i)=p_1$ such that $D_i(p_1,k) = 1$.  This value of $p$ is unique because $D_i(p,k)$ exhibits a unique maximum in $(0,1)$. Moreover,  it is very close to 1 for all $k$ and $i$ and tends to 1 as $i$ tends to infinity, meaning that under-dispersion only occurs when the probability of inserting 0 is very low. 

As observed in Figure \ref{varphi}, the standard deviation of $X$ exhibits a maximum at $p=p_m$. Besides, this maximum exists for any $k$ and $i$ and it is achieved for the value of $p$ such that:
\begin{equation}
\frac{d \, VAR ({\bf X}_i)}{d\,p} = 0
\end{equation}
which can be written as follows:
\begin{equation}\label{polyder}
k^i p^{i-1} \left(i + \sum_{j=1}^{i-1}  (i+j) k^{j-1} p^j  -  i k^i p^i \right) = 0
\end{equation}
which implies that:
\begin{equation}
\sum_{j=1}^{i-1} (i+j) k^{j-1} p^j = i \, ( (k\,p)^i - 1) 
\end{equation}
for any $k$ and for each iteration $i$. Summing the finite sum of the left hand side of this equation yields to:
\begin{equation}
\frac{(1 + i) k \, p - i \, k^2 p^2 - 2\, i k^i p^i + (-1 + 2 \, i) k^{1 + i} p^{1 + i}}{k (-1 + k \,p)^2} = i \, ( (k\,p)^i - 1) 
\end{equation}

This equation cannot be explicitly solved for $p$ for any value of $k$ and $i$. Nonetheless, the roots of the polynomial:
\begin{equation}
r_{k,i}(p) =i + \sum_{j=1}^{i-1}  (i+j) k^{j-1} p^j  -  i \, k^i \, p^i
\end{equation}
can be obtained numerically. For each substitution length, $k=1,2,\ldots$, we are able to fit the roots of polynomial $r_{k,i}(p)$ to the function of $i$:
\begin{equation}
r_k(i) = \frac{1}{1+ \frac{\alpha(k)}{\alpha(k)+(i-1)^{\beta(k)}}}
\end{equation}

The values of the fitting parameters $\alpha(k)$ and $\beta(k)$, as well as the corresponding residual sum of squares, are shown in Table \ref{uniquetable}. As it can be observed, as the $k$-value increases, it seems that  $\alpha(k) \to 0.5$ and $\beta(k) \to 1$, that is:
\begin{equation}
r_k(i) \to r(i) = 1 - \frac{1}{2\,i} 
\end{equation}
This suggests that the roots of the polynomials, i.e., the values of $p$ where the variance of $X$ achieves its maximum, tends to 1 as $i$ tends to infinity. In addition, the variance of $X$ at the maximum tends to infinity, although, by definition,  this variance is null at $p=1$.  We have not a rigorous explanation of the reason behind the almost independence with $k$ of the roots of the polynomial $r_{k,i}(p)$. Likely, it seems to be related to the self-similar properties of the generation process.

\begin{table}[h]
\caption{{Values of the} 
 fitting parameters for some values of the substitution length $k$. As it can be seen, as $k$ increases the values of $\alpha(k) $ and $\beta(k)$ seems to converge to 0.5 and 1, respectively. The third column depicts the Residual Sum of Squares ($RSS$)  of the fitted model.}\label{uniquetable} 

\newcolumntype{C}{>{\raggedright\arraybackslash}X}
\setlength{\tabcolsep}{14mm}{
\begin{tabularx}{\textwidth}{llll}
\toprule
{\textbf{k}} 
 & {\textbf{\textalpha}}   & {\textbf{\textbeta}} & {\textbf{RSS}} \\
\midrule
2 &  0.6256     &   1.0645 &  {$1.78\times 10^{-5}$} 
\\
3  &  0.5832     &  1.0486  & {$ 2.35\times 10^{-5} $} \\
4 &    0.5622    & 1.0383  &  {$ 2.10\times 10^{-5}$} \\
5 &     0.5497   &  1.0316   & {$1.73\times 10^{-5}$} \\
10  &   0.5248    & 1.0167    &  {$6.99\times 10^{-6}$}\\
100 &   0.5025 &  1.0017   &  {$1.05\times 10^{-7}$}\\
\bottomrule
\end{tabularx}}
\end{table}


\section{Entropy}

In this section, we study the mean entropy of a population of sequences generated from the substitution rules defined in the previous section.  As a matter of illustration, we consider first the entropy of the Bernoulli process of probability $p$:
\begin{equation}\label{Hchain0}
H(p) = -p\,log_2(p) - (1-p)\,log_2(1-p)
\end{equation}
This function is symmetric and exhibits a maximum at $p=\frac{1}{2}$. If we now consider the sequences formed by $n$ stochastic variables that follow a Bernoulli process, we can compute the entropy of a sequence of length $n$ with a number of ones $j$ in terms of the frequency $\frac{j}{n}$, that coincides, in average, with $p$:
\begin{equation}\label{Hchain}
H(j) = -\frac{j}{n}\,log_2(\frac{j}{n}) - (1-\frac{j}{n})\,log_2(1-\frac{j}{n})
\end{equation}
Obviously, the sequences with minimum information content, $\frac{j}{n} = p = 0$ and $\frac{j}{n}=p=1$,  have null entropy, i.e.,  $H(0)=H(1) = 0$. In contrast, the maximum entropy is achieved for $j =\frac{n}{2}$ and takes the value: $H(\frac{n}{2}) = 1$.

For a population of Bernoulli sequences of probability $p$, the mean entropy can be calculated applying the binomial distribution to the entropy of the sequences generated with a probability $p$:
\begin{equation}
H^B_n(p) =  {\bf H} \, {\bf P}(p)  = \sum_{k=0}^{n} \,  H(k) \, \binom{n}{k}\, p^k \, (1-p)^{n-k} 
\end{equation}
where ${\bf H}$ is the row vector entropy for a sequence of length $n$ whose components, $H(k)$ are the corresponding entropies Equation (\ref{Hchain}) for all possible combinations of ones. For instance, for $n=1$, this formula reduces to:

\begin{equation}
H^B_1(p)=\binom{2}{0}\,p^0(1-p)^2 \, H(0) +\binom{2}{1}\,p^1(1-p)^1 \, H(1) + \binom{2}{2}\,p^2(1-p)^0 \, H(2) = 2 \, p \, (1-p)
\end{equation}

{where $H(k)$ are the entropies of sequences with 0, 1 and 2 ones: $H(0) = 0, H(1)  =  log_2(2) = 1$ and $H(2) = 0$. Similarly, the mean entropy for a length $n=2$ can be calculated as follows}: 


\begin{equation}
H^B_2(p) = \sum_{k=0}^3 \, \binom{3}{k} \, p^k \, (1-p)^{n-k} \, H_2(k) = -3 \, p \, (1-p) \left( (1-p) \left(\frac{3}{4} \, \log(3) -2 \right) + p \right)
\end{equation}

To compute the mean entropy of a population of sequences randomly generated by substitutions according to the rules defined in Section \ref{sec2}, we have to know the probability of appearance of each of the frequencies of ones in order to weights adequately the sequence entropies. As an example, let us start with the substitution length $k=2$. In this case, the probability distributions of 1 for the $i$-th iteration are given by Equation \eqref{nite}, and this yields the mean entropy: 
\begin{equation}\label{GenH}
H_i(p) = {\bf H}_i \, {\bf P} ({\bf X}_i) 
\end{equation}
Here, the components of the entropy vector ${\bf H}_i$ are given by:
\begin{equation}
H_{i}(j) = - j\,2^{-i} \log_2(j\,2^{-i}) - (1-j\,2^{-i}) \, \log_2 ((1-j\,2^{-i}))
\end{equation}
for a sequence of length $n=2^i$, with a number of ones $j$, such that $0 \leq j \leq 2^i$.

Using the matrix expression for ${\bf P}$, (Equation \eqref{Mexp}), this equation reduces to:

\begin{equation}
H_i(p) =   {\bf H}_i \, M_{(2^i+1)\times(2^{i-1}+1)}\,  M_{(2^{i-1}+1)\times(2^{i-2}+1)} \ldots   M_{5\times3} \, M_{3\times1}  \, {\bf P} ({\bf X}_1)= {\bf H}_{i} \, M_{(2^i+1) \times 1}\, {\bf P} ({\bf X}_1)
\end{equation}

The entropy in the first iteration, i.e., after substituting the initial 1 by applying the rule (\ref{Grule}) with $k=2$ is given by:
\begin{equation}
H_1(p) = \left(0, -\frac{1}{2} \,log_2\left(\frac{1}{2}\right) - \frac{1}{2}\,log_2\left(\frac{1}{2}\right),0\right) \, \begin{pmatrix} \binom{2}{0} \, (1-p)^2 \\ \binom{2}{1} \, p\,(1-p) \\ \binom{2}{2} p^2 \end{pmatrix}  
\end{equation}
which yields:
\begin{equation}
H_1(p)= 2 \,p \, (1-p) 
\end{equation} 

In the next iteration, the entropy can be calculated by the formula:
\begin{equation}
H_2 (p) = {\bf H}_2 \, M_{5\times3}   \,{\bf  P} ({\bf X}_1) 
\end{equation}
where

\begin{equation}\footnotesize
{\bf H}_2 = \left(0, - \frac{1}{4} \,log_2\left(\frac{1}{4}\right) - \frac{3}{4}\,log_2\left(\frac{3}{4}\right), - \frac{1}{2} \,log_2\left(\frac{1}{2}\right) - \frac{1}{2}\,log_2\left(\frac{1}{2}\right), - \frac{3}{4} \,log_2\left(\frac{3}{4}\right) - \frac{1}{4}\,log_2\left(\frac{1}{4}\right), 0\right)
\end{equation}

Multiplying the matrices, we get:
\begin{equation}
H_2(p) = 10 \, p^2 \,(1 - p) \, \left((p^3 - p^2 + \frac{1}{5}\,p + \frac{4}{5}) - \frac{3}{5}\,(p^3 - p^2 + \frac{1}{2})\,log_2(3) \right)
\end{equation}

In the same way, we can calculate the entropy of the 8-digit binary sequences obtained after three substitutions:

{\scriptsize 
\begin{align}
&H_3(p) = 202 \, p^3\,(1 - p) \nonumber \\
&\left[-20/101\,p^3 + 4/101\,p + p^{10} - 3\,p^9 + 285/101\,p^8 + 15/101\,p^7 - 173/101\,p^6 + 91/101\,p^5 + 17/101\,p^4 - 12/101\,p^2 + 12/101 \right. + \nonumber \\
&+  \left(-7/101\,p^{10} + 21/101\,p^9 - 42/101\,p^8 + 35/202\,p^7 + 42/101\,p^6 - 42/101\,p^5 - 7/101\,p^4 + 21/202\,p^3 + 7/101\,p^2 - 7/202 \right)\,log_2(7) - \nonumber \\
&- \frac{35}{101} \, p \left(-27/35\,p^7 + 9/14\,p^6 + 39/35\,p^5 - 33/35\,p^4 - 9/35\,p^3 + 3/35\,p^2 + 9/70\,p + 3/70 \right)\,log_2(3)  + \nonumber \\
&+ \frac{35}{101} \, p^3 \, (1-p)  \left( p^6 - 2\,p^5 + 4/7\,p^4 + 13/14\,p^3 - 1/2\,p^2 - 3/14\,p - 1/14 \right) \,log_2(5)\Big]
\end{align}
}

\begin{figure}[h]
\includegraphics[width=0.7\textwidth]{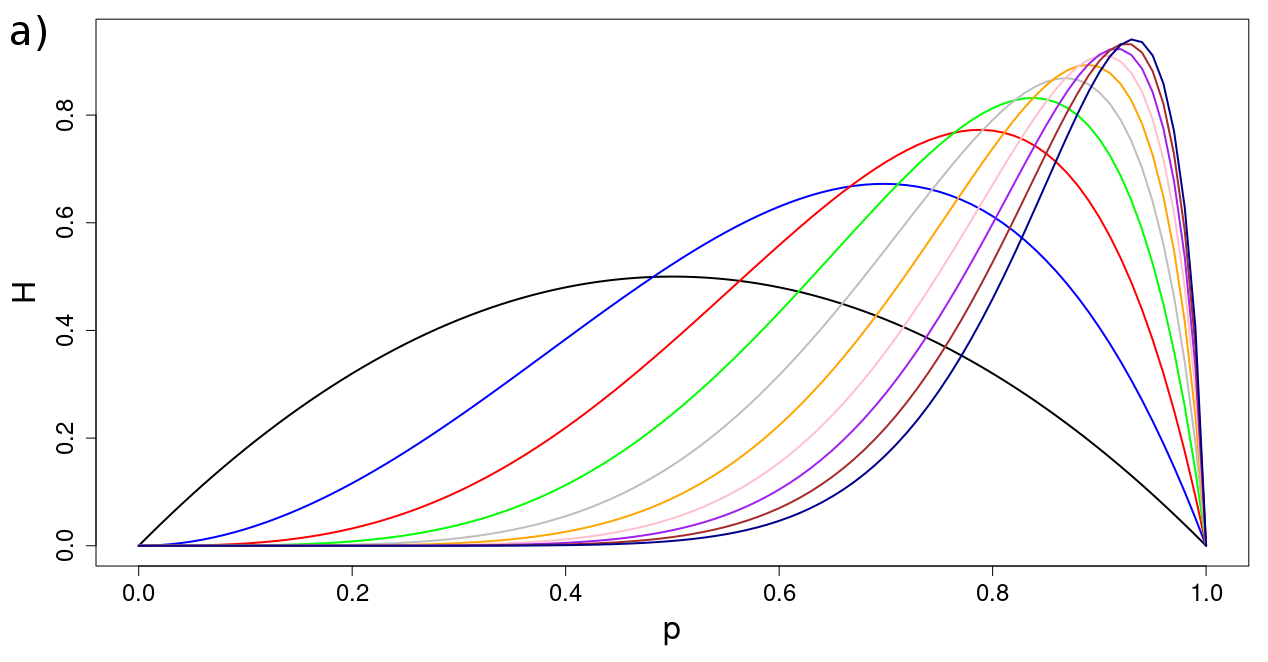}
\includegraphics[width=0.7\textwidth]{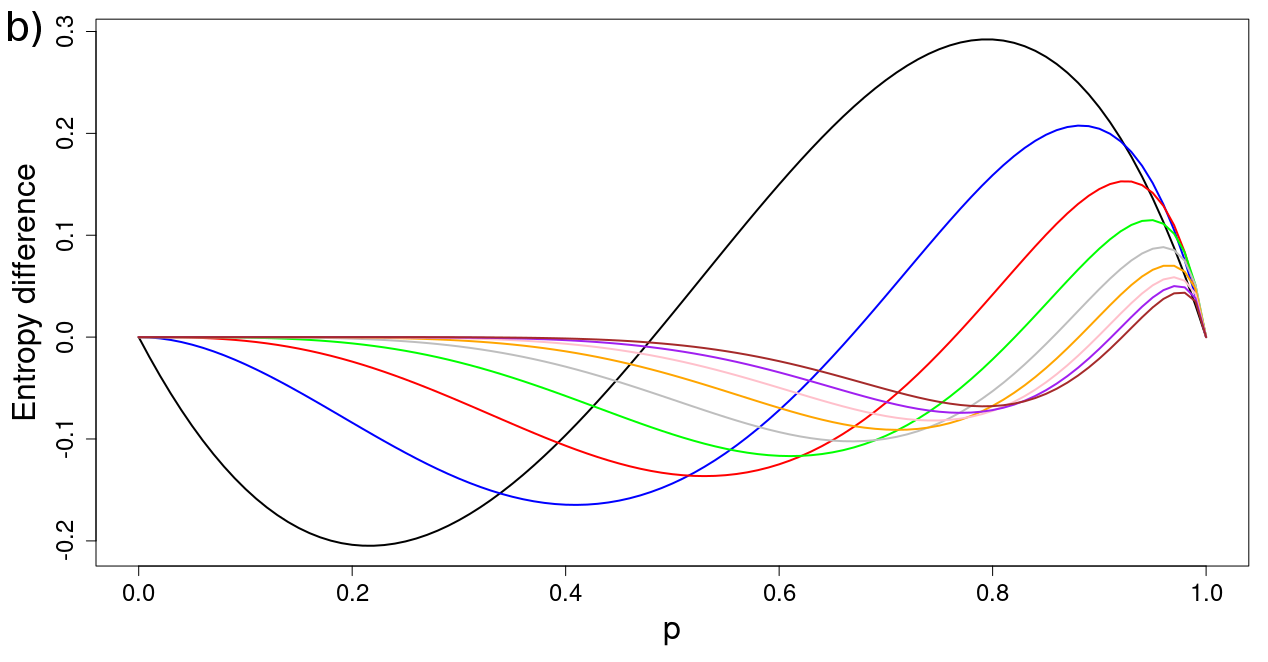}
\caption{({\bf a}) {Entropy as} 
 a function of $p$ for the first 10 iterations. Note the displacement of the maximum to larger values of $p$ as the $i$ increases.  ({\bf b}) Differential entropy $h_i(p) = H_{i+1}-H_i$, for $i =1,2,3,\ldots,9$. As it can be seen, $h_i(p)$ is positive for $p>p_{ch}(i)$. As the iteration increases, the interval of $p$ where $h_i$ is positive shrinks. As a matter of fact, in the limit, it tends to 0. Furthermore, $h_i$ tends to 0 as $i \to \infty$, which means a convergence to a limit entropy for the infinite sequence. The substitution length is $k=2$.}
\label{FigHVAR1}
\end{figure}

In Figure \ref{FigHVAR1}, we depict the entropies calculated analytically for the first three iterations, as well as their corresponding numeric estimations.  Note that, contrary to the Bernoulli sequences, the entropies of these sequences are not symmetric and achieve a maximum at a value of $p$ that depends on the iteration. Furthermore, it can be shown that the entropy "per digit", i.e., $\frac{1}{2^i} \, H_i$, tends to 0 as $i \to \infty$ for all values of $p$.

The differential entropy for successive iterations: 
\begin{equation}
h_i(p) = H_{i+1}(p) - H_i(p) \hspace*{1cm} \text{for} \hspace*{0.5cm} i=1,2,\ldots
\end{equation}
presents an interesting behaviour, as it can be seen in Figure \ref{FigHVAR1}b. There is a critical value that separates two regimes: for $p < p_{ch}(i)$, $h_i(p)$ is negative, whereas for $p >p_{ch}(i)$, $h_i(p)$ is positive. As it can be seen, $p_{ch}(i)$ depends on the iteration: as $i \to \infty$, $p_{ch}(i) \to 1$. As a consequence, the $p$-interval where $h_i(p)$ is positive shrinks to 0. Besides, $h_i(p)$ tends to 0 for all $p$, which means that the entropy converges to the function $H_{\infty}(p)$, for the infinite sequence
.

\begin{figure}[h]
\includegraphics[width=0.8\textwidth]{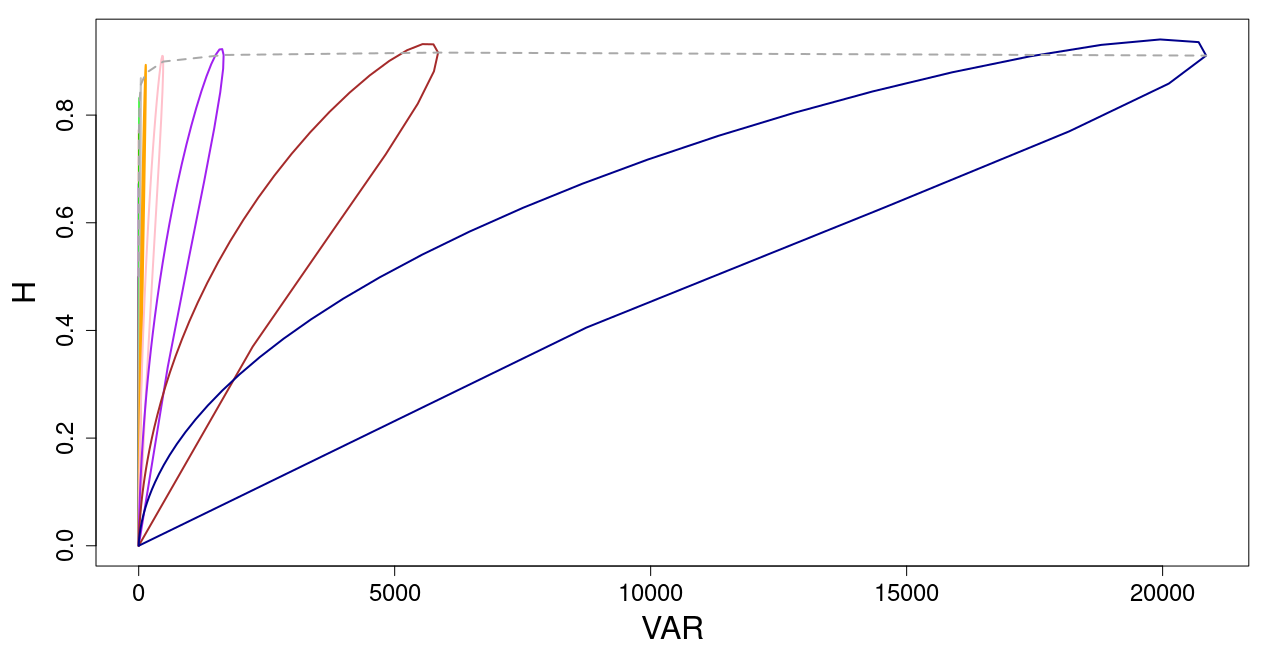}
\caption{{Parametric} 
 curves $H-VAR$ in terms of $p$  for the first 10 iterations. All these curves can be split into two regimes: (i) for $0<p<p_r$, where $H$ has a concave (down) dependence on $VAR$ and (ii) for $p >p_r$, when both variables decreases with $p$ almost linearly.  The value of $p_r$ coincides with the point of the curve farthest from the origin. The dashed line connects these points. Due to the large value of the variance in comparison with the corresponding entropy values, the value of $p_r$ tends to be the value where the maximum of VAR occurs. An interesting consequence of this dependence is that, for a given iteration $i$, we can find two values of $p$ that yield sequences with the same variance but, with different entropies,  having one of them the maximum value. The substitution length is $k=2$.}
\label{FigHVAR2}
\end{figure}

Because both variance and entropy are measures of uncertainty, it could be interesting to find a functional dependence between both (see Figures \ref{FigHVAR1} and \ref{FigHVAR2}). The curves we see, for the case $k=2$, can be parametrized either by $p$, for each iteration $i$, or by $i$, for each $p$. The parametric equations for each $p$ and $i$ are given by Equations \eqref{GVAR} and \eqref{GenH}, for the variance and the entropy, respectively.  As it can be observed, the relationship between $H$ and $VAR$ depends on $p$ for each iteration $i$. It is important to remark that, contrary to what is known for classical probability distributions \cite{Ebrahimi, Mukherjee}, the dependence of $H$ on $VAR$ is not the graph of a function. Instead, the variation of $p \in [0,1]$ brings about a closed curve. Certainly, for $p < p_r$, this functional dependence coincides with that of the classical probability distributions \cite{Mukherjee}. On the contrary, for $p>p_r$,  where both variables decrease with $p$, we find an almost linear dependence between them.  Figure \ref{FigHVAR2} depicts the complete parametric curves for the first iterations obtained analytically.

An interesting consequence of this dependence is that, for the same variance, we can find a value of $p$ that maximizes the entropy. This can be done for any iteration $i >1$ (remember that for $i=1$, $H(p)=VAR(p)$ for all $p \in [0,1]$). Note that the dispersion index, $D$, of these sequences would be different because the mean would change at both values of $p$.  As an example, we find two values of $p$ with a similar variance, in particular, $VAR(0.85) \approx VAR(0.99) \approx 8700$ but, with a mean entropy values $H(0.85) =0.67$ and $H(0.99)=0.41$.  Obviously, we can also find values of $p$ that give rise to sequences with the same entropy that have different variance. For example, a similar mean entropy is achieved for $p=0.79$ and $p=0.99$, whereas the variance of these populations are very distant, $VAR(0.79) \approx 3368$ and $VAR(0.99) \approx 8741$.

\section{Concluding Remarks}

Motivated by the classical percolation process defined by Mandelbrot \cite{Mandelbrot}, we have studied a class of random substitutions of constant length that maps $0 \to \langle 0, 0, \ldots, 0  \rangle$ and $1 \to \langle Y_1, Y_2, \ldots,Y_k \rangle$, being $Y_j$ a Bernoulli process with parameter $p \in [0,1]$.  Contrary to the classical Bernoulli process, this is an asymmetric rule that yields infinite sequences with an average ratio $\#1/\#0$ that tends to 0 for any $0\leq p < 1$. By applying this substitution map, we have randomly generated binary sequences whose length increases with the substitution length $k$ as $k^i$ for $i=0,1,2,\ldots$ for iteration $i$.  Specifically, we have analyzed the combinatorial, statistical and entropic properties of 
this substitution process starting from a sequence formed by the digit 1.

This asymmetry is reflected in  the mass distribution of the number of ones, specifically with regards to its uncertainty.  We have computed the variance of this distribution and the entropy of the ensemble of sequences generated at each iteration as a function of $p$. The relationship between these two magnitudes is depicted in the H-VAR curves of Figure \ref{FigHVAR2}. It is interesting to remark  the two regimes that appear in the curves H-VAR, i.e., concavity down for $p<p_r$ and up for $p >p_r$ that, to the best of our knowledge, have not been shown before \cite{Mukherjee}. This allows to compare sequences generated with different values of $H$ but with the same variance and, on the contrary, sequences that have the same entropy and different variance.  We have left for a forthcoming paper the study of other properties of the generated sequences as, for instance, the distribution of the sizes of the sets of ones and the fractality, measured by the Haussdorff dimension.

The same kind of substitution rules can also be performed in the plane or in the space to generate sets with specific properties. In particular, random substitutions have been applied to study percolation properties of 2D fractal objects \cite{Mandelbrot, Dekking}. Using an iterated function system \cite{Falconer}, random fractals are generated as a function of a probabilistic control parameter, e.g. the probability of any site of being occupied, and their percolation characteristics are quantified in terms of this parameter.


\begin{thebibliography}{99}

\bibitem{Fogg} Fogg,N.P; Berthé, V.;  Ferenczi, S.;   Mauduit, C.; Siegel, A. (Eds.) {\it  Substitutions in Dynamics, Arithmetics and Combinatorics}; Springer: {Berlin/Heidelberg, Germany,} 
2002.


\bibitem{Dekking} Dekking, F.M.; Meester, R.W.J. On the Structure of Mandelbrot's Percolation
Process and Other Random Cantor Sets. \emph{J. Stat. Phys.} \textbf{1990}, {\it 58}, 1109--1126.



\bibitem{Falconer} Falconer, K. {\it Fractal Geometry: Mathematical Foundations and Its Applications}; John Wiley \& Sons Ltd.: New Jersey, NJ, USA,
1990.


\bibitem{Mandelbrot} Mandelbrot, B.B. {\it The  Fractal  Geometry  of  Nature}; W. H. Freeman and Company: {New York, NY, USA}, .
1983.


\bibitem{Barbe} Barbé, A.; Von Haeseler, F.; Skordev, G. Limit sets of restricted random substitutions. \emph{Fractals} \textbf{2006}, \emph{14},  37–47.


\bibitem{Mukherjee} Mukher jee, D.; Ratnaparkhi, M.V. On the functional relationship between entropy and
variance with related applications. \emph{Commun. Stat.-Theory Methods} \textbf{1986}, \emph{15}, 291--311


\bibitem{Ebrahimi} Ebrahimi, N.; Maasoumi, E.;  Soofi, E.S. Ordering univariate distributions by entropy
and variance. \emph{J. Econom. }\textbf{1999}, \emph{90}, 317--336.


\bibitem{Morales} García-Morales, V. Substitution systems and nonextensive statistics. \emph{Phys. A Stat. Mech. Appl.} \textbf{2015}, \emph{440}, 110--117.

\bibitem{Wolfram} Wolfram, S. {\it A New Kind of Science}; Wolfram Media Champaign:{ Illinois, IL, USA, } 
2002.


\bibitem{Shallit} Shallit, J.; Stolfi, J.  Two methods for Generating fractals. \emph{Comput. Graph.} \textbf{1989}, {\it 13}, 185--191.



\bibitem{Glebeta} Gel$\beta$, P.; Schütte, C. Tensor-generated fractals: Using tensor decompositions for creating self-similar patterns. \emph{arXiv} \textbf{2018}, arXiv:1812.00814v1.


\bibitem{Voevudko} Voevudko, A.E. Fractal dimension of the kronecker product based fractals.  \emph{arXiv} \textbf{2018}, arXiv:1803.02766v1.


\bibitem{Xue} Xue, D., Zhu, Y., Zhu, G. X., and Yan, X. Generalized kronecker product and fractals. In \emph{SPIE Digital Library, Proceedings of the Fourth International Conference on Computer-Aided Design and Computer Graphics, {Wuhan, China, 23--25 October 1995};} International Society for Optics and Photonics: {Washington, DC, USA.}


\bibitem{Schroeder} Schroeder, M.R. {\it Fractals, Chaos, Power Laws: Minutes from an Infinite Paradise}; {Courier Corporation: Massachusetts, MA, USA,} 
 1991.




\bibitem{Baake} Baake, M.; Grimm, U.; Penrose, R. {\it  Aperiodic Order}; Cambridge University Press: Cambridge, UK, 2013.


\bibitem{Chayes} Chayes, L. Aspects of the fractal percolation process. Progress in Probability, Vol. 37. In {\it Fractal Geometry and Stochastics}; Springer: {Berlin/Heidelberg, Germany,} 
1995.




\bibitem{oeis} Sloane, N.J.A. Sequence A202019 in The On-Line Encyclopedia of Integer Sequences. Avaliable online: {http://www.oeis.org} {(accessed on February 16th, 2022)} 

\end{thebibliography}
\end{document}